
 The following paper is typed in AMS-Tex Version 2.1.

\input amstex
\magnification =1200
\documentstyle{amsppt}
\refstyle{A}
\widestnumber\key{222}

\NoBlackBoxes

\hsize =16.0truecm
\baselineskip =16.00truept
 \vsize =23.58truecm
\parskip 4pt
\hcorrection{0.17truein}

\leftheadtext{YUN-Gang YE}
\rightheadtext{Complex Conformal Manifolds}
\TagsOnRight

\topmatter
 \title Extremal Rays and Null Geodesics \\
 on  a complex conformal manifold \endtitle
 \author  Yun-Gang Ye\endauthor
 \address
 Department of Mathematics,
  Texas A\&M University, College Station, TX. 77843. 
 Current Address:  Department of Mathematics, Purdue University,
 West Lafayette, IN. 47907, USA.
 \endaddress

\endtopmatter

\def\DD{\nabla}
\def\C{\Bbb C}
\def\D{D}

\def\lr{\longrightarrow}

\def\Q{\Bbb Q}
\def\P{\Bbb P}
\def\v{\varphi}

\def\a{\alpha}
\def\b{\beta}

\def\o{\Omega_X^1}

\def\t{T_X}

\def\ga{\sum g_{\alpha\,ij}dz_{\a}^idz_{\a}^j}

\def\lp{L^{\bot}}
\def\po{\P\left(\o\right)}
\def\pt{\P\left(\t\right)}
\def\op{\Cal O}

\def\pf{\P^1\overset f\to\lr X}
\def\pv{\P^1\overset v\to\lr X}

\def\Ga{\Gamma^i_{\a\,jk}}
\def\intprod{\mathbin{\hbox{\vrule height .5pt width 3.5pt depth 0pt %
      \vrule height 6pt width .5pt depth 0pt}}}

\document

 \heading{\bf $\S$0. Introduction}\endheading

    Let $X$ be a $n$-dimensional complex manifold. A  holomorphic
	conformal structure
	on $X$ is an everywhere {\it non-degenerate} holomorphic section
	$g\in H^0\left(S^2\Omega_X^1(N)\right)$ for some holomorphic line bundle $N$
on $X$.
	Locally, $g$  can be thought of as
	a holomorphic metric  given by $g_{\alpha}=\ga$ on  the
	coordinate chart $U_{\alpha}$ such that $det(g_{\alpha\,ij})$ is everywhere
	non-zero. On the overlap $U_{\alpha}\cap U_{\beta}$, we have
	$g_{\alpha}=f_{\alpha\beta}g_{\beta}$, where $f_{\alpha\beta}\in
	 \op^*_{U_{\alpha}\cap U_{\beta}}$ is an
	invertible holomorphic function on $U_{\alpha}\cap U_{\beta}$. The
	set of holomorphic functions $\{f_{\a\b}\}$'s are the transition
	functions for the line bundle $N$. We call $N$ the {\it conformal
	line bundle} of $X$. It is
	clear that any  complex torus admits a holomorphic conformal structure
	 with trivial conformal line bundle.
	A more interesting example is a smooth hyperquadric
	$\Q^n\subset \P^{n+1}$ (see [5], or $\S$2  for a  description).

		  About ten years
	ago, after classifying all compact complex
	 conformal surfaces and conformal manifolds
	 of any dimensions with  K\"ahler-Einstein metrics,  Kobayashi and Ochiai
	  proposed in [5] the following question:

	\proclaim{Question}  Let $X$ be a compact complex manifold of
	dimension $n$ with
	a holomorphic conformal structure. If $c_1(X)>0$ and $n\geq 3$,
	 is it true that  $X\cong \Q^n$ ?.
	\endproclaim

	 This question was answered positively in [5] assuming the existence
	 of a K\"ahler-Einstein metric on $X$. Their proof was based on
	 Berger's holonomy reduction theorem.  When
	 $n$ is odd, the answer to this question is trivially positive
	  because of index considerations (see (1.1)).

	  This paper grew out of an attempt to understand
	  relationships between extremal rays and null geodesics on a
	  complex conformal manifold. As a consequence, we are able to
	   show that the answer to the above question is positive. In fact,
	    we will prove
	   a stronger result. Precisely, we will show:
	   \proclaim{Theorem 2.1} Let $X$ be a $n$-dimensional
	    complex projective manifold with a holomorphic conformal structure. If
		$n\geq 3$ and $K_X$ is not nef, then  $X\cong \Q^n$.
	   \endproclaim

 Complex  conformal geometry plays an important role in the Penrose's
	   twistor program.  The compactified (and complexified) Minkowski space
	   $\overline {\C M}$
	   is
	    a smooth four-dimensional hyperquadric
		 $\Q^4$ with a natural conformal
		structure.  Once a
		conformal structure (in any dimension) is given, we can define  the notion of
{\it
		 complex null
		geodesics}.  In the case of complexified Minkowski space,
		these correspond to  (complex) light rays.
		In general, null geodesics  are holomorphic curves which are both null (this
means
		their tangent vectors are null vectors)  and geodesics (with
		respect to the Levi-Civita connections of local representatives
		of  the  conformal  metric). The notion of null geodesic is
		well-defined globally on a conformal manifold since two conformally
		equivalent metrics have
		the same null cones and null geodesics (see [2] for a proof).
		For a precise definitions of these notions, we refer
		 readers to [2], or $\S$1.  LeBrun  studied the
		space of all null geodesics on a general complex conformal manifold.
		 Assuming that  the space of null geodesics  is globally convex
		 (this means that any two points on the manifold can be connected
		 by at most one null geodesic), he
		 was able to give a pretty good description of the space. For
		 example, he showed that
		 the space has a natural contact structure.

		  Theorem 2.1
		 has the following interesting  consequence.

	\proclaim{Corollary 1.7} If $X$ is a $n$-dimensional
	($n\geq 3$) complex projective conformal
manifold with an immersed  rational null geodesic, then $X\cong \Q^n$.
\endproclaim

		 The idea of the proof of Theorem 2.1 comes from two seemingly unrelated
		 sources:  Mori's theory of extremal
		 rays
		 and  LeBrun's work on null geodesics. We first  give two
algebraic-geometric  characterizations of null geodesics (see Corollary 1.4
 and Corollary 1.6
 below).  These last two results  enable us to use the local convexity of
the space of null geodesics to conclude that a family
of certain rational curves are locally convex. By {\it local convexity}
of a family of curves
 we mean that there  is at most one curve in the family
 connecting two given points in a neighborhood (in classic topology)
   around any point of $X$.  Local convexity of null geodesics
    is guaranteed by a theorem of Whitehead (see [8]).
 The use of local convexity is an essential part
   of our proof.

  Here is the strategy of the proof.
We
first show that $X$ is a Fano manifold with Picard number one, and
 there exists a divisor in $X$  whose intersecting number with a
 rational curve is one. This implies that the index
 of the manifold is $\text{dim}\, (X)$, which in turn implies the
 Theorem 2.1  immediately by virtue of a theorem of Kobayashi and Ochiai
[4]. The
divisor we will construct is the locus of all the  rational curves
 (or equivalently,  null geodesics) with certain numerical properties
 in $X$ passing through a fixed point.  In the case when $X\cong \Q^n$, this
 divisor is simply the intersection of $X$ with its tangent hyperplane at a
point.

     The paper is organized as follow. There are two sections. In the first
section, we study
some general properties of a complex conformal manifold.  The
main goal of that section is to establish a criterion for
 extremal rational curves to be null geodesics.
The second section  is devoted to the proof of Theorem 2.1.

 All the varieties in this paper are defined over $\C$.  Almost all our
notations are
standard. We will explain them when they occur. But there are several
frequently used
notations which we fell that it is necessary to explain  at this point.
For any coherent sheave ${\Cal F}$, we denote by ${\Cal F}^*$ its dual.
For a line bundle, or an invertible sheave,  $\Cal L$, both $2{\Cal L}$ and
${\Cal L}^2$ mean
${\Cal L}\otimes {\Cal L}$. If
$E$ is a vector bundle, or a  locally free sheave, then we define $\P(E)$ to be
$E^*\setminus \{0\}/\C^*$,  instead of $E\setminus \{0\}/\C^*$. The reason
for this is that it is more conventional for  algebraic geometers.
 $S^mE$ ($m$ is a positive integer) means the $m$-th symmetric tensor
 of the vector bundle $E$.
We  simply denote by $\op(l)$  ($l$ is an integrer) the line bundle
$\op_{\P^1}(l)$
on $\P^1$. If $\Cal F$ is a coherent sheave on $\P^1$, we denote
 ${\Cal F}\otimes\op_{\P^1}(l)$ by ${\Cal F}(l)$.
 We denote by $T_X$  (respectively, $\Omega_X^1$) the  holomorphic
 tangent bundle (respectively, cotangent bundle) of $X$.
All morphisms $f:C\lr X$ from a smooth curve into a  complex manifold $X$ are
assumed to be {\it bimeromorphic}  onto its image.
Throughout this paper we will use the theory of extremal rays freely.
 We refer readers to [1] and [7] for basic materials
 about this subject.

\noindent{\bf Acknowledgments:} The author thanks
Robert McLean and Philip Yasskin for very interesting discussions.

 \heading{\bf $\S$1. Algebraic Characterizations of Null Geodesics}\endheading

 To start with, let $X$ be a complex holomorphic conformal manifold with
a conformal structure $g\in H^0(S^2\Omega_X^1(N))$.
We sometimes call $g$ a (holomorphic) conformal metric on $X$. Let $n$ be the
dimension of $X$ through this paper.
Let $C$ be a smooth complex curve, and $C\overset f\to\lr X$ be a
holomorphic immersion, i.e., the image $f(C)$  have  at worst unramified
singularities.

 We are going to introduce some definitions.
 These definitions are explained  in details in [2], to which we refer
the readers.

Let $v$ be a type (1,0) tangent vector on $X$, $v$ is called {\it null}
if $g(v,v)=0$. Since multiplying by a conformal factor does not change the
 equation $g(v,v)=0$, the notion of nullity is  well-defined globally.
At each point $x\in X$, we have a null tangent cone denoted by  $\Q_x\subset
T_{x,X}$.  It is the set of all null tangent vectors
of $X$ at the point $x$. This was called the {\it sky} at $x$.
It is the affine cone over a smooth hyperquadric.  An immersed holomorphic
curve
$C\overset f\to \lr X$ is called {\it null} if for any point $p\in f(C)$ on
the curve,
and any holomorphic tangent direction $v_p$  at $p$, $g_p(v_p,v_p)=0$.
This definition can be  extended to any  map $f:C\lr X$,
not necessarily unramified. In general, $f(C)$ is called null
if the regular part of $f(C)$ is null. This clearly agrees with
the definition  for immersed curves.

On each local coordinate chart $U_{\a}$,  let
 $g_{\a}=\ga$ be a local representative of the conformal metric $g$.
 We can associated to $g_{\a}$ a Levi-Civita
connection  $\DD_{\a}$ on $U_{\a}$, very much like the case  in the real
Riemannian geometry.  Namely, we have the Christoffel symbols
$$\Ga=\frac{1}{2}\sum_m
g^{im}_{\a}\left(\partial_{j}g_{\a\,km}+\partial_{k}g_{\a\,jm}
-\partial_{m}g_{\a\,jk}\right)$$
where $\partial_k=\partial/\partial z_{\a}^k$, and
$[g_{\a}^{ij}]=[g_{\a\,ij}]^{-1}$
 as matrices.  An immersed holomorphic curve $C$ in $U_{\a}$ is a (complex)
 {\it geodesic} if there exists a holomorphic parameter $\xi$ on $C$ such that:
 $\dfrac{d}{d\xi}\intprod\DD_{\a}\left(\dfrac{d}{d\xi}\right)=0$, i.e.,
$$
\dfrac{d^2z_{\a}^i}{d\xi^2}+\sum_{j,k}\Ga\dfrac{dz_{\a}^j}{d\xi}\dfrac{dz_{\a}^k}{d\xi}=0$$
for all $i=1,\dots, n$.

But these connections on various coordinate charts
 do not match  on the overlaps since the conformal metric $g$ is not globally
defined.
However, since   they  differ only by some conformal
factors on the overlaps, the notion of  null geodesic is well-defined globally
(see [2] for a proof).
For a morphism $f: C\lr X$,
$f(C)$ is called a {\it null geodesic} if it is both null and a geodesic on
each coordinate cover $U_{\a}$.

 Let $Y=\pt\overset\text{def.}\to = \Omega^1_X\setminus\{0\}/\C^*$.
   Using the conformal
structure $g$, we obtain a well-defined divisor of $Y$, which we denote
 by $S$.  Precisely,
$$S=\bigg\{(p,\omega_p)\, \bigg\vert\,p\in X\,, \omega_p\in \Omega^1_{p,X}\,,
 g_p(\omega_p,\omega_p)=0\bigg\}$$
where $\Omega^1_{p,X}$ is the holomorphic cotangent space of $X$ at the point
$p\in X$, and $g_p$ is any representative of the conformal metric $g$ at
the point $p$.
 Therefore $S$ is a smooth quadric bundle over $X$.
Let $Y\overset \pi\to\lr X$ be the natural projection, and $S_x=\pi^{-1}(x)
\cap S$ for a point $x\in X$.  Then it is
clear that for any $x\in X$,  $S_x\cong \Q^{n-2}$ and the sky $\Q_x$
is an affine cone over $S_x$.

Denote by $L$ the tautological line bundle $\op_Y(1)$ on $Y$.
\proclaim{Lemma 1.1}  $\op_Y(S)\cong 2L\otimes\pi^{*}N^{*}$.
\endproclaim
\demo{Proof}  Note that the conformal structure induces a bundle
isomorphism:
$T_X\overset g\to\cong\Omega_X^1(N) $.
Let  $g^{-1}\in H^0(S^2T_X\otimes
N^*)$ be
the inverse of $g\in H^0(S^2\Omega_X^1(N))$. In view of the isomorphism
 $$H^0(L^2\otimes \pi^*N^{*})\cong H^0(S^2T_X\otimes N^*)$$
  $g^{-1}$ induces  a
section $s_g\in H^0(L^2\otimes N^{*})$. It is straightforward to
check that the divisor $S$ constructed above is exactly the
vanishing locus of $s_g$, i.e., $S=(s_g)_0$. This implies
the lemma immediately.

\qed
\enddemo

 The following fact about a complex conformal manifold will be used
 frequently in the rest of the paper.
 $$ 2K_X=-nN \tag 1.1$$
  where $K_X=\wedge^n\Omega_X^1$ is the canonical bundle of $X$. (1.1)  can be
proved by taking the determinant  of the isomorphism
 $T_X\cong \Omega_X^1(N)$.

Next we will give some  algebraic characterizations of
null-geodesics on a complex conformal manifold. First let us observe that
 there is a natural contact structure on $Y=\pt$.  The associated
 contact line bundle is $\op_Y(-1)$, which is the dual of
 $L$. Therefore we have the
following short exact sequence:
$$0\lr \lp \lr T_Y\overset \phi\to\lr L\lr 0 \tag 1.2$$
The contact form is given locally by $\theta=\sum_{i=1}^{n}\xi_idz^i$,
where
 $z=(z^1,\dots,z^n)$ is a local coordinate on $X$, and
 $[\xi_1,\cdots,\xi_n]$
 a  projective coordinate on the fiber of $\pi: Y\lr X$. More precisely,
on the chart, say, where $\xi_1\neq 0$,  the contact form
is given by
$\theta_1=dz^1+\sum_{i=2}^{n}\left(\dfrac{\xi_i}{\xi_1}\right)dz^i$. Therefore
the bundle $\lp$
is spanned by all tangent vectors perpendicular to the contact
form $\theta$. The homomorphism
 $\phi$ is defined via contraction by the contact form.
 The  one form $\theta$, in fact,  lives naturally on the whole
cotangent bundle $\Omega_X^1$ and $d\theta=\sum_{i=1}^{n}d\xi_i\wedge dz^{i}$
is a natural symplectic form on $\Omega_X^1$.

Choose a
coordinate cover $\{U_{\a}\}$ of $Y$ such that $\theta_{\a}$ is the local
 contact form. Let $\omega_{\a}=d\theta_{\a}\vert_{\lp}$. Then
$\{\omega_{\a}\}$'s define a {\it conformal symplectic } structure on
the bundle $\lp$.  That is, $\{\omega_{\a}\}$'s may not be glued to
give a well-defined symplectic form on $\lp$, but  on the overlap
 $U_{\a}\cap U_{\b}$, we have:
  $\omega_{\a}=h_{\a\b}\omega_{\b}$,
where $h_{\a\b}$ is the transition function for the contact line
bundle $\op_Y(-1)$ on $Y$.  The equation
 $\omega_{\a}=h_{\a\b}\omega_{\b}$ is
obtained by differentiating the equation $\theta_{\a}=h_{\a\b}\theta_{\b}$
 (since $\theta_{\a}|_{\lp}=0$ for each $\a$).

In terms of these local coordinates, we can represent
 the conformal structure locally by an invertible symmetric $n\times n$ matrix
 $\left(g^{ij}\right)$. And the divisor $S\subset Y$ can be locally defined by
the
 equation
   $\{((z^i),[\xi_i])\vert \sum_{1\leq i,j\leq n}g^{ij}(z)\xi_i\xi_j=0\}$.
If we restrict the bundle map $\phi$ in (1.2) to
 $T_S\subset T_Y\vert_S$, we get a homomorphism of bundles:
$T_S\overset \phi\vert_S\to\lr L\vert_S$. It is easy to check that  $\phi
\vert_S$ is
 surjective using local coordinates.
Let $T_S^0$ be $\text{ker}\, (\phi\vert_S)$. Then $T_S^0$ is a well-defined
subbundle of $\lp\vert_S$. The rank of $T_S^0$ is $2n-3$ since  $\lp$ is of
 rank $2n-2$. Denote $\omega$ by the
collection of $\{\omega_{\a}\}$'s.  Restricting
$\omega$ to $T^0_S$,    we get a bundle map  locally on each coordinate chart:
$\omega(S): T_S^0\lr T^{0*}_S$. Even though $\omega(S)$ is not well-defined
 globally, the kernel of $\omega(S)$ is. Let $N^0_S\overset \text{def.}\to =
\text{ker}\,(\omega(S))$.

 We claim that $N^0_S$ is a  line subbundle of
$T^0_S\subset \lp|_S$. This can be shown as follow. First of
 all,  $\omega(S)$ has to degenerate
 since  rank of $T_S^0$ is  ($2n-3$), which is odd. Secondly, the rank of
$\text{ker}(\omega(S))
 \leq \text{codim}(S)=1$ since  each $\omega_{\a}$ is non-degenerate. Hence
 $\text{ker}(\omega(S))$ is a line subbundle of $T^0_S$.

  Following LeBrun,
 holomorphic integral curves
of  $N^0_S$ are called {\it phase-space trajectories} of $S$.  The importance
of
this notion is illustrated by the following result of LeBrun (see p.213
in [2] for a proof).

\proclaim{Proposition 1.2} [LeBrun] Null geodesics of $X$ are precisely
  the  phase-space trajectories of $S$.
 \endproclaim

 The line bundle $N_S^0\subset T_S^0$ can be identified with $N^*_{S/Y}$,
 i.e.,
 $$N^0_S=N_{S/Y}^{*} \tag 1.3$$
 where $N^*_{S/Y}$ is the  co-normal bundle
of $S$ in $Y$.  We can show (1.3) as follow. The contact sequence (1.2)
induces the following   short
exact sequence:
$$ 0\lr T^0_S\lr \lp\vert_S\lr N_{S/Y}\lr 0\tag 1.4$$
 Using the conformal symplectic structure $\{\omega_{\a}\}$ on $\lp$, we get
 (1.3) immediately.

Now we are going to translate the above result of LeBrun
into the language of algebraic geometry.
Let $C\overset f\to\lr  X$ be a holomorphic  immersion (i.e., unramified)
from a {\it smooth} algebraic curve $C$.
 For any point $p\in f(C)$ and any tangent direction $v_p$ of $f(C)$ at $p$,
 we get a point $(p,v_p)\in \po$. In this way we construct a  lift of
  the map $f$
to $\po$. Denote this lift by $f_1: C\lr \po$.  Let $C_1=f_1(C)\subset \po$.
Then $C_1$ is a smooth algebraic curve since $f$ is an immersion.
  Since $X$ is conformal,  $T_X\otimes N^{-1}\cong \o$.
   This implies that   $\po\cong \pt=Y$.  Now composing this last
   isomorphism with the lift
   $f_1:C\lr \po$,  we get a
lift $f_2: C\lr Y$. Then it is clear that $f_2$ maps $C$ isomorphically
 onto its image.  We still denote  by $C$  the image $f_2(C)\subset Y$ when
 there is no danger of confusion. Therefore the resulting morphism $C\overset
 f \to\lr X$
 is nothing but the projection $C\overset \pi|_C\to\lr X$, i.e,
 $f=\pi|_C$.   It is clear
from the construction  that
 $$L\vert_{C}\cong \Omega_{C}^1\otimes \pi^*N\vert_{C}\tag 1.5$$
 since  $L\otimes \pi^*N^*$ is the tautological
line bundle of $\po$ under the isomorphism $\po\cong \pt$.

  There is an isomorphism similar to (1.5) that holds for a general
  (not necessarily immersed)
  morphism $C\overset f\to \lr X$. Let $B$ be the ramification
  divisor for the morphism $f$. It is the zero divisor of the differential
  of $f$. Then the following is true:
  $$f_2^*L\cong \Omega_C^1(-B)\otimes f^*N \tag 1.5a$$
  where $f_2: C\lr Y$ is the lift of $f: C\lr X$ to $Y$ as constructed
  above.     We can proof (1.5a) as follow. First note that $L\otimes \pi^*N^*$
is
the tautological line bundle of $\P(\Omega_X^1)$. Secondly,
the lift of $C\overset f\to\lr X$ to $\P(\Omega_X^1)$ is defined by the
quotient line bundle
  $\Omega_C^1(-B)$ of $f^*\Omega_X^1$. Therefore
$f_2^*L\otimes f_2^*\pi^*N^*\cong \Omega_C^1(-B)$, i.e.,
$f_2^*L\otimes f^*N^*\cong \Omega_C^1(-B)$ since $f=\pi\circ f_2$. This
implies (1.5a) immediately.
  In general, the lift  $f_2$ does not map $C$ isomorphically
  onto its image due to the presence of ramifications.

  Recall that a curve $C\subset Y=\pt$ is called a {\it contact} curve
  if it is an integral curve for the contact distribution  $\lp$ on $Y$.
  Another way to say is that all the local contact forms are zero when
  they restrict to $C$.

 \proclaim{Lemma 1.3} Suppose that  $C\overset f\to\lr X$ is a holomorphic
 immersion from a smooth algebraic curve.  If we denote still by
 $C$  the lift of $C$ to $Y$ as constructed above, then $C$ is a contact curve
of $Y$  if and only if $f(C)\subset X$ is
  null.
  \endproclaim
 \demo{Proof} Since this is a local question, we may as well assume
 that $f$ is an embedding. Choose local coordinates as above.
 Let $[\xi^1,\dots,\xi^n]$ be the corresponding fiber coordinate for $\po$.
  Then
 $$ \xi^i=\sum_{j=1}^{n}g^{ij}\xi_j,\,\,\,\,\,\,\,\,\,\text{for}\,\,\,\,\,
 \,\,\,1\leq i\leq n\tag 1.6$$
  Let $t$ be a local holomorphic parameter on $C$.
  Then for $1\leq i\leq n$
   the following holds on $C$:
  $$ \xi^i=\frac{dz^i}{dt}\tag 1.7$$
   Since $\theta=\sum_{i=1}^n\xi_idz^i$,
  $$
  \aligned
  \theta\vert_C &=\left(\sum_{i=1}^{n}\xi_i\xi^i\right)dt\\
         &=g((\xi^i),(\xi^i))dt
  \endaligned \tag 1.8$$
  where $g(\cdot,\cdot)$ is a local representative of the conformal structure.
   At this point, the lemma is obvious.

 \qed
 \enddemo

 By definition,  $C$ being contact is equivalent to the fact that
  $T_C$ is a line subbundle of $\lp|_C$.  This last lemma implies that
  $T_C$ is a line subbundle of $\lp|_C$
 if and only if $f(C)$ is a null curve. Note also that the $f(C)\subset X$ is
{\it null} if and only if
 the lift  of $C$ to $Y$  is contained in $S$.
Now Proposition 1.2 implies the following corollary immediately.

\proclaim{Corollary 1.4} Let $C\overset f\to\lr X$ be a
holomorphic immersion from a smooth curve $C$ into $X$ such that
its image is null. Then $f(C)$ is a
null geodesic if and only if the two line subbundles
 $T_C\hookrightarrow \lp\vert_C$ and $N_S^0|_C\hookrightarrow \lp\vert_C$
are identical.
\endproclaim

   Nullity is usually easy to
  check as we see in the next lemma.
  \proclaim{Lemma 1.5}   Let $C\overset f\to\lr X$ be an
holomorphic map (not necessarily an immersion) from a smooth curve $C$ into
$X$. Then $f(C)$ is null if
 $h^0(S^2\Omega_{C}^1\otimes f^*N)=0$.
 \endproclaim
 \demo{Proof} Let $g\in H^0(S^2\Omega_X^1(N))$ be the given conformal
 structure. Then $f(C)$ is null if and only if $f^*g=0$ on $C$.
 Since $f^*g\in H^0(S^2\Omega_{C}^1\otimes f^*N)$,
 the lemma is obviously true.

 \qed
 \enddemo

 The following corollary will play an essential role in the proof of
 Theorem 2.1.
\proclaim{Corollary 1.6} Let $f: C\lr X$ be a null-immersion, i.e.,
 its image is null. If $f(C)$ is a null geodesic,
then $T_C\cong f^*N$. Conversely, if
$T_C\cong f^*N$ and $h^0\left(\lp|_C\otimes \Omega_C^1\right)=1$,
then $f(C)$ is a null geodesic.
\endproclaim

\demo{Proof}  Let us first prove that
$$N^0_S\vert_C\cong
 T_C^{\otimes 2}\otimes \pi^*N^*\vert_C\tag  1.9$$
   By Lemma 1.1, we have:
$\op_Y(S)\cong L^2\otimes \pi^{*}N^{*}$. Since
$N_{S/Y}\cong\op_Y(S)\vert_S$, (1.5)
implies that $N_{S/Y}|_C\cong\left(\Omega_C^1\right)^{\otimes 2}\otimes \pi^*N
\vert_C$. Now (1.9) easily follows from (1.3).

 If $C$ is a null geodesic, then  $T_C\cong N_S^0|_C$ by
 Corollary 1.4.
 By (1.9),   $N^0_S\vert_C\cong
 T_C^{\otimes 2}\otimes \pi*N^*$. Therefore
    $T_C\cong \pi^*N|_C$. However, $f=\pi\vert_C$.
Therefore $T_C\cong f^*N$.

If $T_C\cong f^*N$, then $T_C\cong N^0_S\vert_C$ by (1.9).
If moreover $h^0(\lp\vert_C\otimes \Omega_C^1)=1$, then any line subbundle
 of
$\lp\vert_C$ that is isomorphic  to
$T_C$ has to be indentical to $T_C$ as
a line subbundle of $\lp\vert_C$.       By Corollary 1.4, $C$ is a null
geodesic of $X$.
Hence we are done.

\qed
\enddemo

 If $f:C\lr X$ is not an immersion, there is still an isomorphism similar
 to (1.9). As before, let  $B$ be the ramification locus of $f$. Then
 $$ f_2^*N^0_S\cong T_C^{\otimes 2}(2B)\otimes f^*N \tag 1.9a$$
 where $f_2: C\lr Y$ is the lift of $f$ to $Y$.  This can be easily proved
using (1.5a).

\proclaim{Corollary 1.7} If $X$ is a $n$-dimensional
 ($n\geq 3$) complex projective conformal
manifold with an immersed  rational null geodesic, then $X\cong \Q^n$.

\endproclaim

\demo{Proof}  Let $\pf$ be an immersed rational  null geodesic.   Then
Corollary 1.6
implies that $T_{\P^1}\cong f^*N$.
Therefore $N\cdot f(\P^1)=2$. By (1.1), $2K_X=-nN$. Therefore  \linebreak
$K_X\cdot f(\P^1)<0$,
i.e., $K_X$ is not nef.  Theorem 2.1 below implies the corollary immediately.

\qed
\enddemo

\heading{\bf $\S$2.  Projective Conformal Manifolds with
 Non-Nef Canonical Bundles}\endheading

 Throughout this section, we let $X$ be a $n$-dimensional
complex   projective conformal manifold with  $K_X$ being not nef.
Recall that we say a line bundle $\Cal L$ is {\it nef} if
$\Cal L\cdot C\geq 0$ for any effective curve $C\subset X$. By a celebrated
theorem of Mori ([6] \& [7]),  there is an extremal ray
$R$ on $X$.   From now on, we fix an extremal ray $R$ on $X$.

We define the {\it length} of the extremal ray $R$ (denoted by $l(R)$)
as:
$$l(R)=\bigg\{-K_X\cdot \Gamma\,\bigg\vert\, [\Gamma]\in R\, \,
		   \text{and}\,\, \Gamma\,\,\text{is}\,\, \text{a}\,\,
		   \text{rational}\,\, \text{curve}.\bigg\}$$
		   where $[\Gamma]$ means the numerical equivalence class of
		   the rational curve $\Gamma$.

 By a theorem of Mori [6],
  $l(R)\leq n+1$. We
 say  a rational curve $\Gamma$ is an
  {\it extremal rational  curve} of $R$ if  $[\Gamma]\in R$ and
   $-K_X\cdot \Gamma=l(R)$.

\noindent{\bf Example.} \hskip 3pt Let $X=\Q^n\subset \P^{n+1}$  be
a smooth hyperquadric and $R$ be the  extremal
ray generated by a straight line $\ell\subset X$. Then  $l(R)=n$
and all extremal
rational curves of this extremal ray are straight lines in $X$. There is a
 natural conformal structure on $X$, which can be described as follow.
 Assume, without loss of generality,  that
 $$X=\bigg\{[z_0,\dots\, z_{n+1}]\bigg\vert \sum_{i=0}^{n+1}z_i^2=0\bigg\}
 $$
 Then the symmetric two form $g=\sum_{i=0}^{n+1}dz_idz_i$  on $\C^{n+2}$
 descends to
 a conformal structure on $X$ with the conformal bundle $N\cong\op_X(2)$.
 It is easy to show that this is the {\it only} conformal structure on $\Q^n$.

Now let us state the main theorem in this section.
\proclaim{Theorem 2.1} Let $X$ be a $n$-dimensional
complex projective manifold with
a holomorphic conformal structure. If $n\geq 3$ and $K_X$ is not
nef, then $X\cong
\Q^n$.
\endproclaim

  The proof of this theorem requires several lemmas.  Here is the
 outline of our proof. We first show that $X$ is a Fano manifold
   with Picard number one
  and admits  an extremal ray $R$ such that  $l(R)=n$. We then show that there
  is a divisor $\D$ and a rational curve $C$ on $X$ such that
  $\D\cdot C=1$. This implies that the index of $X$ is $n$.  A theorem of
   Kobayashi and Ochiai (see the next paragrah for a statement of the theorem)
implies that
  $X$ is necessarily a smooth hyperquadric. The holomorphic
  conformal structure
  will play an essential role in our proof. A crucial idea here is
  to show that general  extremal rational  curves are necessarily
   null geodesics (see Lemma 2.7 below). This enables us to show the local
convexity of
   extremal rational curves.

  Let us recall a theorem of Kobayashi-Ochiai [4]. If
  $X$ is a Fano manifold (i.e., $-K_X$ is ample) of dimension $n$, then
   we define the index
  of $X$ (denoted by $\text{index}(X)$) to be the largest positive
   integer that
  divides $K_X$ in the Picard group  $Pic(X)$.
  Then the theorem of
  Kobayashi-Ochiai asserts that $\text{index}(X)\leq n+1$, and
  if $\text{index}(X)=n+1$, then $X\cong \P^n$; if $\text{index}(X)=n$, then
  $X\cong \Q^n$.

   Let $r=0,1$.  Fix $r$ point $\{t_1,\dots,t_r\}$ on $\P^1$.
   Of course, when $r=0$, the set $\{t_1,\dots,t_r\}$ is empty.
Let $\P^1\overset f\to\lr X$ be a morphism such that it is birational
onto its image.
  Consider the space  $Hom(\P^1,X;t_i,f(t_i))$ of all morphisms from $\P^1$ to
$X$
  fixing the given $r$ point. Let $U_r$ be the irreducible
  component of $Hom(\P^1,X;t_i,f(t_i))$ that contains the point $[f]$.
  Here we denote by $[f]$ the point corresponding to the morphism
  $f:\P^1\lr X$.

   Then we have
  natural morphism:
  $$U_r\times \P^1 \overset \Xi\to\lr  X,\,\,\,\,\,\,
\,\,\,\,\,\,\,\Xi([v],t)=v(t)$$
  for any  $[v]\in U_r$ and $t\in \P^1$.
  The following  lemma is a slight variation of Corollary 1.3 of [3].
  \proclaim{Lemma 2.2} Let $d$ be the dimension of the image of $\Xi$.
  Then for a general $[v]\in U_r$, $v^*T_X$ has the following decomposition:
  $$v^*T_X\cong \oplus^{n}_{i=1}\op(a_i)$$
  such that $a_1\geq\cdots\geq a_d\geq r$.

  \endproclaim
 \demo{Proof} Let $G_r$ be the subgroup of
  $Aut(\P^1)$ fixing the given $r$ point. Then $G_r$ acts transitively on
  $\P^1$. Therefore  the image of
   $\Xi\vert_{U_r\times \{t\}}: U_r\times\{t\}\lr X$ also has dimension $d$ for
   any $t\in \P^1$. By the Generic Smoothness
   Theorem,   for a general $[v]\in U$, the differential of
   $\Xi\vert_{U_r\times \{t\}}$ at $[v]\times \{t\}$ has rank $d$. That is,
    the
   following homomorphism
   $$H^0(v^*T_X(-r))\lr v^*T_X(-r)\vert_{\{t\}}\tag 2.1$$
   has rank $d$. This is enough to imply the lemma.

   \qed
   \enddemo

 Let us fix a morphism $\pf$ such that $[f(\P^1)]$ belongs to the
 given extremal ray $R$.  We denote
  by $E(f)$ the image $\Xi(U_0\times \P^1)$, and $F_x(f)$ by
  image $\Xi(U_1\times \P^1)$, where $x=f(t_1)$ for some point $t_1\in \P^1$.
Note that
  $E(f)$ is simply  the  union of the images of all  the deformations of $\pf$.
  $F_x(f)$ is  the union of the images of all the deformations of $\pf$ fixing
  the point $p_1$. It is clear to see that
  $E(f)$ is contained in the exceptional locus $E_R$ of the contraction
  morphism $\v_R$ of $R$, and $F_x(f)$ is contained in a fiber of
   $\v_R\vert_{E_R}$. The following
   lemma was proved implicitly in the proof of Theorem (1.1) of [10].
   \proclaim{Lemma 2.3} If $-\text{deg}f^*K_X=l(R)$, then:
   $$\text{dim}E(f)+\text{dim}F_x(f)\geq n+l(R)-1\tag 2.2$$
   \endproclaim
   Since the proof of the lemma is not very difficult, we will give a proof
below.
   \demo{Proof} There is a natural morphism $\a: U_0\lr Chow^{l(R)}(X)$, the
   Chow variety of 1-cycles with intersection number $l(R)$ with $-K_X$, such
   that $\a([f])=f(\P^1)$ for any $[f]\in U_0$.
   Let $T_0$ be the image of $\a(U_0)$. Consider the following
   morphism:
   $$\P^1\times U_0 \overset \Psi \to\lr T_0\times X,\,\,\,\,\,\,\,\,\,\,\,
\,\,
    \Psi(t,[f])=(\a([f]),\Xi(t,[f]))$$
	Let $M_0\subset T_0\times X$ be the image of $\Psi$. Then we have the
following diagram:
	$$
	\CD
	M_0 @>p>> X\\
	@VqVV \\
	T_0
	\endCD
	$$
	where $p$ and $q$ are two projections.

	 First note that
	$E=p(M_0)$.  Let  $T_0(x)=q(p^{-1}(x))$. Then
	$F_x(f)=p\left(q^{-1}T_0(x)\right)$.
	If we
	 denote still  by $p$ its restriction to $q^{-1}T_0(x)$, then
	$p: q^{-1}T_0(x)\lr F_x(f)$ is finite. Otherwise, there would be a
	 non-trivial
	deformation of $h(\P^1)$ for some $[h]\in U_0$
	 fixing  two different points on $h(\P^1)$. By the
	breaking-up technique (see [6]), there must be a rational curve
	$\Gamma\subset X$ such that $[\Gamma]\in R$ and $-K_X\cdot \Gamma
	< -K_X\cdot h(\P^1)=l(R)$. This is a contradiction by the definition
	of the length of an extremal ray.  Therefore 	$p: q^{-1}T_0(x)\lr F_x(f)$ is
finite.

	Hence $\text{dim}F_x(f)=\text{dim}q^{-1}T_0(x)=\text{dim}T_0(x)+1$.
  Since $q$ maps $p^{-1}(x)$ isomorphically onto $T_0(x)$, we have:
  $$
  \aligned
  \text{dim}E &\geq \text{dim}M_0-\text{dim}T_0(x)\\
               &= \text{dim}M_0-\text{dim}F_x(f)+1
  \endaligned
  $$
  However
  $$
  \aligned
  \text{dim}M_0 &=\text{dim}U_0+1-\text{dim}Aut(\P^1)\\
                &\geq \chi(f^*T_X)+1-3\\
				&= n+l(R)-2
  \endaligned
  $$
   Here  we use Proposition 3 in  [6] to conclude that
  $\text{dim}U_0\geq \chi(f^*T_X)$.
  Combining these inequalities, we obtain the lemma immediately.

   \qed
   \enddemo

   \proclaim{Lemma 2.4}  Let $X$ be a $n$-dimensional  projective conformal
manifold with
   non-nef canonical bundle. Then there is an extremal ray $R$ on $X$
   such that $l(R)=n$.
 \endproclaim
 \demo{Proof}  By the Cone Theorem (see [7]), there exists an extremal ray
 $R$ on $X$.  Since $2K_X=-nN$ (by equation (1.1))  and $l(R)\leq n+1$, $l(R)$
is either $n$, or
  $\dfrac{n}{2}$ ($n$ is necessarily even in this case). If $l(R)=n$, then
  we are done.  Assume  that $n=2k$ and $l(R)=k$.  We will
  derive a contradiction in this case.

  Let $\pf$ be a morphism such that $f^*K_X=\op(-k)$.  Then
  $N\cdot f(\P^1)=1$
 by (1.1). For simplicity, let $E$ be $E(f)$.
   Without loss of generality, we can assume that $x=f(0)$ is a smooth
   point of $f(\P^1)$.  Denote $F_x(f)$ by $F_x$. Then Lemma
  2.3 implies that
  $$\text{dim}E\geq 3k-1-\text{dim}F_x\tag 2.3$$

  \noindent{\it Claim:} \hskip 4pt $\text{dim}F_x=k$.
  \demo{Proof of the claim} Let $d=\text{dim}F_x$.
  Denote by $U_1\ni [f]$ an irreducible component of $Hom(\P^1, X; 0,x)$.
  For a general
  point $[u]\in U_1$, let $u^*T_X\cong \oplus_{i=1}^{2k}\op(a_i)$ such that
  $a_1\geq\dots\geq a_{2k}$. By Lemma 2.2, we have
   $$a_1\geq\dots\geq a_d\geq 1\tag 2.4$$
  However the conformal structure  $g$ induces an isomorphism
  $T_X\overset g \to\cong \Omega_X^1(N)$. This last isomorphism implies:
  $$\oplus_{i=1}^{2k}\op(a_i)\cong \oplus_{i=1}^{2k}
  \op(1-a_i)\tag 2.5$$
  Therefore $a_i+a_{2k-i+1}=1$ for $1\leq i\leq k$. Now (2.4) easily implies
  that $d\leq k$.  However,  $d\geq k-1$ by (2.3) (since $\text{dim}E\leq 2k$),
and the equality holds
   if and only if $E=X$.

Suppose that $d\neq k$.  Then $d=k-1$ and $E=X$.
   By Lemma 2.2,  for a general point $[v]\in U_0$, if
   $v^*T_X=\oplus_{i=1}^{2k}\op(a_i)$, then $a_i\geq 0$ for all $i$.
   In particular, $a_{2k}\geq 0$.  However
   $a_{2k}=1-a_1$ by (2.5) and $a_1\geq 2$ since $\P^1\overset v\to\lr X$ is
not a constant
   morphism. Therefore $a_{2k}<0$. This is a contradiction. Hence
   $d=k$ and $E$ is a divisor by (2.3). Hence the claim is proved.

  \enddemo

  The proof of the claim implies that $E$ is  a divisor.
    Again let  $\P^1\overset v\to\lr X$  be a general point in $U_0$.
   Suppose
    that $v^*T_X\cong \oplus_{i=1}^{2k}\op(a_i)$ such that
	$a_1\geq\dots\geq a_{2k}$.  Then  by Lemma 2.2,
	at least  the first $2k-1$ (if not all) $a_i$'s are
	non-negative. In view of  (2.5), $a_2=\dots=a_{k}=1$, $a_{k+1}=\dots
	=a_{2k-2}=0$ and $a_{2k}=1-a_1$.
	 Hence:
	 $$v^*T_X\cong \op(a_1)\oplus^{k-1}\op(1)\oplus^{k-1}\op\oplus \op(1-a_1)\tag
	  2.6$$
   Since   $\P^1\overset v\to\lr X$ is  not constant,
  $a_1\geq 2$.

	Let $B$ be the ramification divisor of the morphism $v$. Then
	(2.6) implies that $\text{deg}(B)=a_1-2$, i.e. $\op_{\P^1}(B)\cong
	\op(a_1-2)$. This can be proved by observing that the image of the
tangent bundle homomorphism $T_{\P^1}\cong \op(2)\overset dv_*\to\lr v^*T_X$
lies entirely
in the first component $\op(a_1)$.

 Following the same notations as in the previous
	section, we
 denote by\linebreak  $v_2: \P^1\lr Y=\pt$ the lift of $v$ to $Y$.  Therefore
$v=\pi\circ v_2$.
	 Let $L=\op_Y(1)$ be the tautological line bundle of $Y$.
	 By (1.5a), $v_2^*L\cong \op(1-a_1)$ since
	 $v^*N\cong \op(1)$.
	  Consider the relative Euler sequence:
	  $$0\lr L^*\lr \pi^*\Omega_X^1\lr T_{Y/X}\lr 0\tag 2.7$$
	  where $T_{Y/X}$ is the relative tangent bundle for $Y=\pt \overset
	  \pi\to \lr X$. Now (2.6) and (2.7) imply that:
	  $$v_2^*T_{Y/X}\cong \op(-a_1)\oplus^{k-1}\op(-1)\oplus^{k-1}\op \tag 2.8$$

	  On the one hand, since $a_1\geq 2$,  we have
	  $$h^0\left(v^*T_X\otimes \op(1-2a_1)\right)=
	  h^0\left(v_2^*T_{Y/X}\otimes \op(1-2a_1)\right)=0 \tag 2.9 $$
	  by
	  (2.6) and (2.8).
	However $T_Y$ fits into the following exact sequence
	   $$0\lr T_{Y/X}\lr T_Y\lr \pi^*T_X\lr 0$$
	   Therefore     (2.9) implies that
	   $h^0\left(v_2^*T_Y\otimes \op(1-2a_1)\right)=0$.

	   On the other hand,   Lemma 1.5 implies that   $v(\P^1)$ is a null
	    curve since
\linebreak 	$S^2\Omega_{\P^1}\otimes v^*N\cong \op(1-2a_1)$ and $1-2a_1<0$
    (remember $a_1\geq 2$).
	 Hence
	 $C$ is contained in $S$, the quadric bundle over $X$ as defined in the
previous
	 section. As we see in the previous section that  $N_S^0$ is a line
	  subbundle of $T_S\subset T_Y\vert_S$.
	  Therefore $v_2^*N^0_S$ is a line subbundle of $v_2^*T_Y$.
	  Note that $v_2^*N_S^0$ is isomorphic to $\op(2a_1-1)$  by
	   $(1.9a)$. Hence
	   $h^0\left(v_2^*T_Y\otimes \op(1-2a_1) \right)>0$.   This contadicts what we
just
	   prove in the previous paragraph. Hence we are done, i.e., $l(R)$ must
	   be $n$.

 \qed
 \enddemo

  By Lemma 2.4, the extremal $R$ on $X$  has length $n=\text{dim}
X$, i.e., $l(R)=n$. Fix an
  arbitrary
  point $x\in X$.
  Let us fix a  morphism $\pf$  such that $f(0)=x$ is a smooth
  point of $f(\P^1)$, and $f(\P^1)$ is an extremal rational curve of
  $R$, i.e., $K_X\cdot f(\P^1)=-n$.  Lemma 2.3 implies
  that the union of the
 images of all the deformations of $f$ covers the whole space $X$.
  Therefore
   we can assume that $f^*T_X$ is semi-ample by Lemma 2.2.

   Let $V\ni [f]$ be an irreducible component of $Hom(\P^1,X;\{0,x\})$,
   where $[f]$ is the point of $V$ corresponding to the morphism
   $\pf$.
   Let $G$ be the group of automorphisms of $\P^1$ fixing the origin.
   Then $G$ acts  on $V$  by: $g\circ [v](t)\overset\text{def} \to=
   v(g^{-1}t)$ for any $[v]\in V$ and $t\in \P^1$.

	Consider the natural morphism:
	$$\Xi: \P^1\times V\lr X,\,\,\,\,\,\,\,\,\,\,\, \, \Xi(t, [v])=v(t)
	$$
	 Let $D$ be the {\it closure} of $\Xi(\P^1\times V)$.  By Lemma 2.3,
	  it is clear  that $D$ is either a divisor, or $D=X$. We claim that $\D$
	 is necessarily a divisor.   If $\D=X$, then Lemma 2.2 implies that
	 for a general point $[v]\in V$, $a_i\geq 1$ for all $i$'s if
	 $v^*T_X\cong \oplus_{i=1}^{n}\op(a_i)$. Since $a_1\geq 2$, this would implies
	 that $l(R)=-K_X\cdot v(\P^1)\geq n+1$, which contradicts the fact
	 that $l(R)=n$.  Therefore $\D$ is a divisor.

 We will show that $D$ is the
	 divisor we are looking for, i.e., $\D\cdot f(\P^1)=1$.  First, let
	 us prove the following.

   \proclaim{Lemma 2.5} There exists a  neighborhood (in the  classic topology)
    $W$ of $x\in X$ such that $\D\cap W$ is biholomorphic to a neighborhood
	 of the vertex of  the cone
   over a smooth hyperquadric $\Q^{n-2}$ with $x$ corresponding to
   the vertex.
   \endproclaim

 \demo{Proof} Fix a representative, which is denoted by $g_x$,  of the
conformal metric
 $g$ on $X$ locally around $x$.   As it was pointed out in [2] that, by a
theorem of Whitehead [8], we can choose an analytic normal coordinate
neighborhood
 $W \subset_{\text{open}}\C^n$  around the given point  $x$ such that complex
 geodesics (in the metric $g_x$)  are all affine lines in $\C^n$. Let
$z=(z^1\,\dots,z^n)$
  be a     local coordinate. We further assume that
  $x$ is the origin in $B$.
  Let  $Q=\bigg\{z\in W\,\bigg\vert\,g_0(z,z)=0\bigg\}$, where
  $g_0=g_x(0)$.  By Lemma 2.7 below,  $v(\P^1)\subset X$ is a
  null  geodesic through $x$ for a general
point $[v]\in V$. Hence $v(\P^1)\cap W$ is contained in $Q$
   for a general  $[v]\in V$.  By definition, $D$ is the closure  (in Zariski
   topology)
   of  $\bigcup_{[v]\in V}v(\P^1)$.
  Therefore $\D\cap W$ is contained in $Q$, too. But both $Q$ and $\D\cap W$
  have dimensions $n-1$, and $Q$ is an irreducible (since $n\geq 3$) analytic
variety in $B$.
  Therefore $\D\cap W=Q$. This
implies  the lemma immediately.
    Hence we are done.

 \qed
 \enddemo

 \proclaim{Lemma 2.6} Let $X$ be as in Theorem 2.1, then $X$ is a
 Fano manifold with  Picard number one, i.e., $\rho(X)=1$
 \endproclaim

 \demo{Proof} By Lemma 2.4 there exists an extremal ray $R$ on $X$ such
 that $l(R)=n$. In this case, it is easy to show (see
  Proposition (2.4)
 in [9] for a proof) that
 $X$ is either a Fano manifold with Picard number one, or the
 associated contraction map
   $\v_R: X\lr Z$  sends $X$ to a smooth curve $Z$. Moreover
  in the later case,  a general
  fiber  $F$ is a Fano manifold with Picard number one and having
  an extremal ray $R_F$ with $l(R_F)=n=\text{dim}F+1$.

 Now assume to the contrary that the later case occurs.  We will
 derive a contradiction. Let $F$ be a general fiber of the
 fibration $\v_R: X\lr Z$.  Hence $F$ is smooth. Let $\pf$ be a non-constant
 morphism such that $f^*K_X\cong \op(-n)$ and $[f(\P^1)]\in R$. Furthermore we
assume that
 $f(0)=x\in F$. Let $\D$  be the divisor defined before, i.e., the closure
 of
 $\Xi(\P^1\times V)$, where
  $V\ni [f]$ is an irreducible component of $Hom(\P^1,X;\{0,x\})$.   On the one
hand,
  by Lemma 2.5 above,
 $\D$ is locally a quadric cone near the point $x$. Hence $\D$ is
  {\it singular}
 at $x$ since $n\geq 3$. On the other hand, it is easy to see that  images of
  all the deformations of $f$ fixing $x$  are necessarily contained in $F$
since $\v_R$ is a contraction morphism.
 Hence $\D=F$ since they have the same dimensions. However  $F$  is smooth at
$x$. This is a contradiction. Hence
 $X$ must be a Fano manifold with $\rho(X)=1$. We are done.

\qed
\enddemo

  \proclaim{Lemma 2.7} Let $[v]\in V$ be a general point. Then
  \roster
  \item $v^*T_X\cong \op(2)\oplus^{n-2}\op(1)\oplus\op$. In particular,
  $v$ is unramified.
  \item
  $v(\P^1)\subset X$ is a null geodesic.
  \endroster
 \endproclaim

 \demo{Proof}
  Let $\pv$ be a general closed point in $V$.
 Suppose that $v^*T_X\cong \oplus^{n}_{i=1}\op(a_i)$ with
 $a_1\geq\dots\geq a_{n}$. Since rational curves from $V$ cover $X$,
$a_i\geq 0$ for all $1\leq i\leq n$ by Lemma 2.2.  However, the isomorphism
$T_X\cong \Omega_X^1(N)$ and (1.1) imply  that:
$$a_i+a_{n-i+1}=N\cdot v(\P^1)=2 \tag 2.10$$
for $1\leq i\leq n$.

 Since $v:\P^1\lr X$ is non-constant,  $T_{\P^1}\cong \op(2)$ is a subsheave
 of $v^*T_X$. Hence $a_1\geq 2$. Now there are only  two possibilities,
 either $a_2<2$, hence $a_i<2$  for all $2\leq i\leq n$, or
 $a_2\geq 2$. In the former case, we have
 $a_1=2$, $a_2=\dots=a_{n-1}=1$ and $a_{n}=0$
 since $\text{deg}\, v^*K_X=-n$. In the later case,
 $a_1\geq a_2\geq 2$. Let us
 rule out
the later case.

Suppose that  $a_1\geq a_2\geq 2$. By (2.10), $a_1=a_2=2$ since
all $a_i$'s are non-negative. Hence there is an
$\op(2)$ piece in $v^*T_X$ which is different from the    image of
  $T_{\P^1}
\cong \op(2)$ in $v^*T_X$.
 We choose  a global section  $s$ of that $\op(2)$ piece such that $s$ vanishes
  at two distinct points  $t_1$ and $t_2$ of $\P^1$. Moreover we can arrange
   so that both $v(t_1)$ and $v(t_2)$ are smooth points on $v(\P^1)$. Since
$h^1(v^*T_X)=0$,  $V$ is smooth at $[v]$ by Proposition 3 in [6]. Hence
   the section $s$ generates an actual  non-trivial deformation
   of $f(\P^1)$ in $X$ fixing two distinct points, namely, $v(t_1)$ and
$v(t_2)$. By the breaking-up
   technique (see [6]),  there exists a rational
   curve $\Gamma$ such that $[\Gamma]\in R$ and
   $-K_X\cdot \Gamma<-K_X\cdot v(\P^1)=n$.  This contradicts the fact that the
length
   of the extremal ray $R$ is $n$. Hence the later case can not happen.
   The first part of the lemma is thus proved.

 Let $C_v=v(\P^1)$ be the image of $v$ and  $C$ be the lift of
 $C_v$ to $Y=\pt$ as constructed in $\S$1.  Since $N\cdot C_v=2$,
 $S^2\Omega_{\P^1}^1\otimes v^*N\cong \op(-2)$, $C_v$
 is a null curve by Lemma 1.5.

 \noindent{\it Claim:} We have:
  $$ \lp|_C\cong \op(-2)\oplus^{n-2}\op(-1)\oplus\op(2)\oplus^{n-2}
 \op(1) \tag 2.11$$
\demo{Proof of the claim} Consider the relative Euler sequence:
 $$ 0\lr L^{*}\lr \pi^*\Omega_X^1\lr T_{Y/X}\lr 0\tag 2.12$$
  where $T_{Y/X}$ is the relative tangent bundle for
 $\pi: Y=\pt\lr X$, and $L=\op_Y(1)$ is the tautological
 line bundle of $Y$.  Note that $L\vert_C\cong\op$ by (1.5). Now (2.12) implies
that:
 $$T_{Y/X}\vert_C\cong \op(-2)\oplus^{n-2}\op(-1)\tag 2.13$$
 Hence
 $$T_Y|_C\cong \op(-2)\oplus^{n-2}\op(-1)\oplus\op(2)\oplus^{n-2}\op(1)
 \oplus\op  \tag 2.14$$
 since $T_Y$ fits in the following exact sequence:
 $$0\lr T_{Y/X}\lr T_Y\lr \pi^*T_X\lr 0$$
 In view of sequence (1.2), we have:
 $$\lp|_C\cong \op(-2)\oplus^{n-2}\op(-1)\oplus\op(2)\oplus^{n-2}
 \op(1) \tag 2.15$$
  This proves the claim.
\enddemo
 The above claim implies that $h^0\left(\lp\vert_C\otimes \Omega_C^1\right)=1$.
Since $T_{\P^1}\cong\op(2)\cong v^*N$,
  Corollary  1.6 implies the second part of  lemma immediately.

 \qed
 \enddemo

    Let $V_0$ be a Zariski open subset of $V$ such that for any
   $[v]\in V_0$, $v(\P^1)$ is smooth at $x=v(0)$ and Lemma 2.7 holds for $[v]$.
    It is clear that $V_0$ is non-empty.  By Lemma 2.7,
	$h^1(v^*T_X(-1))=0$ for any $[v]\in V_0$. Hence $V_0$ is contained
	in the smooth part of $V$, i.e., $V_0$ is smooth. It is important
	to note that $V_0$ is {\it locally convex} around $x$ by
	Lemma 2.7 due to the local convexity of null geodesics. The local
	convexity of null geodesics is trivial since we can always find
	a normal neighborhood around $x$ (see the proof of Lemma 2.5).

   Consider a natural morphism $\a: V_0:  \lr Chow^{n}(X)$ such that
   $\a([v])=v(\P^1)$ as an  1-cycle.  Let
   $Y_0$ be the normalization of $\a(V_0)$ in $\C(V_0)^G$, where $\C(V_0)^G$
    is the filed of rational functions that are invariant under $G$.
 $\C(V_0)^G$	 is a finite extension of $\C(\a(V_0))$ (see p.602 in [6]).  Then
   the induced morphism $V_0\overset \gamma\to \lr Y_0$ is a geometric
   quotient for the  $G$-action on $V_0$
   by Lemma 9 in [6]. In particular, $V_0$ is a principal
   $G$-bundle
   over $Y_0$, hence $Y_0$ is smooth since $V_0$ is. Note that $Y_0$ may not be
complete since $V_0$ is only
   an open subset of $V$. But this is enough for the proof of Theorem 2.1.

 Now let us prove  Theorem 2.1.

 \demo{\bf Proof of Theorem 2.1}
  Note that the group $G$ acts on $\P^1\times V$ by
  $g(t, [v])=(g(t),g\circ [v])$ for any $g\in G$. Consider the $G$-invariant
morphism
  $$F: \P^1\times V_0\lr Y_0\times X\,\,\,\,\,\,
\,\,\,\,\,\,F(t,[v])=(\gamma([v]), v(t)) $$
  It is clear  that the morphism  defined before
  $\Xi: \P^1\times V\lr X$ is related to  $F$ by
  $\Xi|_{\P^1\times V_0}=p_2\circ F$, where
$p_2$ is the projection to the second factor.

Following Mori [6], let us define $Z_0=\text{Spec}_{Y_0\times X}\left(
\left(F_*\op_{\P^1\times V_0}\right)^G \right)$.
Then  $Z_0$  is the geometric quotient of $\P^1\times V_0$ under the
$G$-action, and
  is a $\P^1$-bundle (in Zariski topology)
over $Y_0$  with a section $S_0=\{([v], 0)\vert [v]\in Y_0\}$. In particular,
  $Z_0$ is  also smooth since $Y_0$ is.  The  above morphism $F$  induces
   a morphism $\Pi: Z_0\lr Y_0\times X$ since $F$ is $G$-invariant.
   Let $\psi_i=p_i\circ \Pi$ ($i=1,2$), where $p_1$ and $p_2$ are the
   two projections for $Y_0\times X$. Then $\psi_2: Z_0\lr \D\subset X$
   is  dominant and the section $S_0$ is contracted to the point
    $x$.    Since
 Lemma 2.7 holds for any $[v]\in V_0$,
  we can use the
 same argument in [6] to conclude that $\psi_2$ is etale
 away from $S_0$.  Namely, the proof goes as follow.
It is enough to show that $V_0\cong \{t\}\times V_0 \overset
p_2\circ F|_{\{t\}\times V_0} \to\lr X$ is
smooth for any $t\in \P^1\setminus \{0\}$.  It is clear that
$p_2\circ F|_{\{t\}\times V_0}([v])=v(t)$.  Therefore
any fiber of $p_2\circ F|_{\{t\}\times V_0}$ is of the form
$Hom(\P^1,X;v|_{\{t,0\}})\cap V_0$ (for some $[v]\in V_0$), the space
of morphisms from $\P^1$ to $X$ fixing the two distinct points $t$ and $0$.
 On the one hand,  $\text{dim}\, Hom(\P^1,X;v|_{\{t,0\}})
\cap V_0\geq 1$
since the the subgroup of $Aut(\P^1)$ fixing the  two points
$t$ and $0$ is one-dimensional.
On the other hand, the Zariski tangent space
of $Hom(\P^1,X;v|_{\{t,0\}})\cap V_0$   at any point $[u]$  is isomorphic
$H^0(u^*T_X(-2))$, which is exactly one-dimensional
since the morphism $u$ satisfies Lemma 2.7  by the definition of
$V_0$. Therefore $Hom(\P^1,X;v|_{\{t,0\}})\cap V_0$ must be smooth
 and  one-dimensional for any $[v]\in V_0$. Hence $V_0\cong \{t\}\times V_0
\overset
p_2\circ F|_{\{t\}\times V_0} \to\lr X$ is
smooth for any $t\in \P^1\setminus \{0\}$.  Therefore, $\psi_2$ is
etale away from $S_0$.

  We claim that $\psi_2: Z_0\lr \psi_2(Z_0)\subset \D$ is in fact isomorphic
  away from $S_0$.  To prove this, we have to make essential use of
  the local convexity of null geodesics. Choose a normal
  coordinate neighborhood $W$  (in classic topology) around $x$ (see the proof
of
  Lemma 2.5).  Therefore, for any $[v]\in V_0$,
$v(\P^1)\cap W$ is an affine line through $x$ by Lemma 2.7. Let
$W_0=\psi_2^{-1}B$. Then $W_0$ is a neighborhood
  (in classic topology) of the section $S_0$.  Pick a point $[v_0]\in Y_0$
  and a point $t_0\in \P^1$ such that $(t_0,[v_0])\in W_0\setminus S_0$
 (that is, $v_0(t_0)\neq x$) and
  $\psi_2(t_0,[v_0])=v_0(t_0)$ is a smooth point on $v_0(\P^1)$. It is clear
that
  $\psi_2^{-1}\psi_2(t_0,[v_0])$ consists of
   only one point by the local convexity of null geodesics.
   Therefore
  $\psi_2:  Z_0\lr \psi_2(Z_0)\subset \D$ is isomorphic away from $S_0$ since
it is etale away from
   $S_0$.

 Hence
   $v(\P^1)$ is smooth away from $x$ for any $[v]\in Y_0$. However $v(\P^1)$
   is already assumed to be smooth at $x$. Hence $v(\P^1)$ is everywhere
   smooth for any $[v]\in V_0$.
   We also have that for any two different points $[v_1]$ and $[v_2]$ in $V_0$,
   $v_1(\P^1)$ and $v_2(\P^1)$ intersect only at $x$. Since two
   different null-geodesics can not tangent to each other,
   $v_1(\P^1)$ and $v_2(\P^1)$ have different tangent directions at $x$.

  Let $\sigma: \widetilde X\lr X$ be the blowing-up of $X$ at $x$. Let $E$ be
the
  exceptional locus of $\sigma$.  By what we prove in the previous
   paragraph, the morphism $\psi_2: Z_0\lr X$ can be lifted to
  $\widetilde \psi_2: Z_0\lr \widetilde X$. Moreover, $\widetilde \psi_2$ maps
$Z_0$
  isomorphically
  onto its image.
    It is clear that $\widetilde \psi_2(Z_0)$ is a Zariski open subset of
	  the proper
  transform $\widetilde \D$ of $\D$ under the blowing-up.

   Now choose an arbitrary  point  $[v]\in V_0$ and let $C=v(\P^1)$. Let
$\widetilde v: \P^1\lr
   \widetilde X$ be the lift of $[v]$ to $\widetilde X$, and
   denote $\widetilde v(\P^1)$ by  $\widetilde C$. Then
   both $\widetilde C$ and $C$ are smooth and   $\widetilde C$ is the proper
transform
   of the curve $C$.
    Since $\widetilde \D$ is isomorphic to $Z_0$ in a neighborhood  of
$\widetilde C$, we have:
	$$ N_{\widetilde C/\widetilde\D} \cong \oplus^{n-2}\op \tag 2.16$$
     It is easy to see that the following sequence is exact:
	 $$0\lr T_{\widetilde X}\overset d\sigma_*\to\lr \sigma^*T_X\lr
\oplus^{n-1}\op_E(-E)
	 \lr 0\tag 2.17$$
	 Now (2.17) implies that
	 $$ N_{\widetilde  C/\widetilde X} \cong \sigma^*N_{C/X}(-1) \tag 2.18$$
	 because $E\cdot \widetilde C=1$ ($C$ is smooth !).
	 Since $[v]\in V_0$,  $v^*T_X\cong \op(2)\oplus^{n-2}\op(1)\oplus\op$
 by the definition of $V_0$. Therefore (2.18) implies that
	 $$N_{\widetilde  C/\widetilde X} \cong \oplus^{n-2}\op\oplus \op(-1) \tag
2.19$$
	Now let us consider the following exact sequence:
   $$ 0\lr N_{\widetilde C/\widetilde \D}\lr N_{\widetilde  C/\widetilde X}\lr
N_{\widetilde \D/\widetilde X}\vert
   _{\widetilde C} \lr 0 \tag 2.20$$
   In view of (2.16) and (2.19), (2.20) implies that
   $N_{\widetilde \D/\widetilde X}\vert_{\widetilde C}\cong \op(-1)$. Hence
  $$\widetilde \D\cdot \widetilde C=-1 \tag 2.21$$
  Since $\D$ is locally isomorphic to a quadric cone by Lemma 2.5, its
multiplicity
  at $x$ is exactly two. Hence:
  $$  \sigma^*\D=\widetilde\D+2E \tag 2.22$$
  Since  $C$ is smooth at $x$, $E\cdot \widetilde C=1$.
   Now (2.21) and (2.22) imply  that
  $\sigma^*\D\cdot \widetilde C=1$, hence $\D\cdot C=1$ by Projection
  Formula.
  Since $\rho(X)=1$ by Lemma 2.4, $\D$ must be the positive generator for
  $Pic(X)$.
  Let $K_X=m\D$. Since $K_X\cdot C=-n$, $m=n$. Hence $K_X=-n\D$, i.e,
  $\text{index}(X)=n$ since $\text{index}(X)\leq n+1$. Now by a theorem of
Kobayashi-Ochiai,
  $X\cong \Q^{n}$.  Hence we complete the proof.

 \qed
 \enddemo


\Refs
\ref \key {\bf 1}
\by  H. Clemens, J. Koll\'ar, S. Mori
\paper Higher Dimensional Complex Geometry
\yr 1988
\jour Ast\'erisque,  Soc. Math. France,
 \vol 166
\endref

\ref\key {\bf 2}\by C. LeBrun \paper Spaces of complex null geodesics
in complex Riemannian geometry
\yr 1983\pages 209-231\jour Trans. Amer. Math. Soc. \vol 278
\endref


\ref\key {\bf 3} \by J. Koll\'ar, Y. Miyaoka, S. Mori
\paper Rationally connected varieties
\jour preprint
\endref




\ref\key {\bf 4} \by S. Kobayashi, T. Ochiai
\paper Characterizations of complex projective spaces and hyperquadrics
\yr 1972
\pages 31-47\jour J. Math. Kyoto Univ. \vol 13
\endref

\ref \key {\bf 5} \by  S. Kobayashi, T. Ochiai
\paper Holomorphic structures modeled after hyperquadrics
\yr 1982\vol 34\pages 587-629\jour Tohoku Math. J.
\endref

\ref \key {\bf 6}\by S. Mori
\paper Projective manifolds with ample tangent bundles
\yr 1979\vol 110\pages 593-606
\jour Ann. Math.
\endref

\ref \key {\bf 7} \by S. Mori
 \paper  Threefolds whose canonical bundles are not
numerically effective
\yr 1982 \pages 133-176
\vol 116\jour Ann. Math
\endref


 \ref \key {\bf 8} \by  J.H.C. Whitehead
 \paper Convex regions in the geometry of paths
 \jour Quart. J. Math. \yr 1932
 \vol 3 \pages 33-42
 \endref

\ref \key {\bf 9} \by  J. Wi\'sniewski
\paper Length of extremal rays and generalized adjunction
\pages 409-427\yr 1987\jour Math. Z. \vol 200
\endref

\ref\key {\bf 10} \by J.  Wi\'sniewski
\paper On contractions of extremal rays on Fano manifolds
\yr 1991\vol 417\pages 141-157
\jour J. reine. angew. Math.
\endref


\endRefs
\enddocument